\documentclass[12pt,leqno]{article}
\usepackage{amsfonts}
\usepackage{amsmath, amsthm, amsfonts, amssymb, color}
\usepackage{mathrsfs}
\usepackage{color}
\usepackage[notcite,notref]{showkeys}

\theoremstyle{definition}

\newcommand{\scr}[1]{\mathscr #1}
\definecolor{wco}{rgb}{0.5,0.2,0.3}

\numberwithin{equation}{section} \theoremstyle{remark}

\newcommand{\ua}{\uparrow}

\title{{\bf Nonlinear Fokker--Planck equations for Probability Measures on Path Space and Path-Distribution Dependent SDEs}\footnote{Supported in
 part by  NNSFC (11771326,11431014).} }
\author{
{\bf    Xing Huang$^{a)}$, Michael R\"ockner$^{b),c)}$, Feng-Yu Wang$^{a),d)}$  }\\
\footnotesize{$^{a)}$ Center for Applied Mathematics, Tianjin University, Tianjin 300072, China }\\
\footnotesize{$^{b)}$ Department of Mathematics, Bielefeld
University, D-33501 Bielefeld, Germany}\\
\footnotesize{$^{c)}$ Academy of Mathematics and Systems Science, Chinese Academy of Sciences, Beijing, 100190, P.R.China}\\
 \footnotesize{ $^{d)}$ Department of Mathematics,
Swansea University, Singleton Park, SA2 8PP, United Kingdom}\\
\footnotesize{  wangfy@tju.edu.cn, F.-Y.Wang@swansea.ac.uk}}
\begin{document}
\allowdisplaybreaks
\def\R{\mathbb R}  \def\ff{\frac} \def\ss{\sqrt} \def\B{\mathbf
B} \def\W{\mathbb W}
\def\N{\mathbb N} \def\kk{\kappa} \def\m{{\bf m}}
\def\ee{\varepsilon}\def\ddd{D^*}
\def\dd{\delta} \def\DD{\Delta} \def\vv{\varepsilon} \def\rr{\rho}
\def\<{\langle} \def\>{\rangle} \def\GG{\Gamma} \def\gg{\gamma}
  \def\nn{\nabla} \def\pp{\partial} \def\E{\mathbb E}
\def\d{\text{\rm{d}}} \def\bb{\beta} \def\aa{\alpha} \def\D{\scr D}
  \def\si{\sigma} \def\ess{\text{\rm{ess}}}
\def\beg{\begin} \def\beq{\begin{equation}}  \def\F{\scr F}
\def\Ric{\text{\rm{Ric}}} \def\Hess{\text{\rm{Hess}}}
\def\e{\text{\rm{e}}} \def\ua{\underline a} \def\OO{\Omega}  \def\oo{\omega}
 \def\tt{\tilde} \def\Ric{\text{\rm{Ric}}}
\def\cut{\text{\rm{cut}}} \def\P{\mathbb P} \def\ifn{I_n(f^{\bigotimes n})}
\def\C{\scr C}      \def\aaa{\mathbf{r}}     \def\r{r}
\def\gap{\text{\rm{gap}}} \def\prr{\pi_{{\bf m},\varrho}}  \def\r{\mathbf r}
\def\Z{\mathbb Z} \def\vrr{\varrho} \def\ll{\lambda}
\def\L{\scr L}\def\Tt{\tt} \def\TT{\tt}\def\II{\mathbb I}
\def\i{{\rm in}}\def\Sect{{\rm Sect}}  \def\H{\mathbb H}
\def\M{\scr M}\def\Q{\mathbb Q} \def\texto{\text{o}} \def\LL{\Lambda}
\def\Rank{{\rm Rank}} \def\B{\scr B} \def\i{{\rm i}} \def\HR{\hat{\R}^d}
\def\to{\rightarrow}\def\l{\ell}\def\iint{\int}
\def\EE{\scr E}\def\Cut{{\rm Cut}}
\def\A{\scr A} \def\Lip{{\rm Lip}}\def\kk{\kappa}
\def\BB{\scr B}\def\Ent{{\rm Ent}}\def\L{\scr L}

\maketitle

\begin{abstract} By investigating path-distribution dependent stochastic differential equations, the following type of nonlinear Fokker--Planck equations for probability measures $(\mu_t)_{t \geq 0}$ on the path space
$\C:=C([-r_0,0];\R^d),$ is analyzed:
$$\pp_t \mu(t)=L_{t,\mu_t}^*\mu_t,\ \ t\ge 0,$$ where
$\mu(t)$ is the image  of $\mu_t$ under the projection  $\C\ni\xi\mapsto \xi(0)\in\R^d$, and
$$L_{t,\mu}(\xi):= \ff 1 2\sum_{i,j=1}^d a_{ij}(t,\xi,\mu)\ff{\pp^2} {\pp_{\xi(0)_i} \pp_{\xi(0)_j }}  +\sum_{i=1}^d b_i(t,\xi,\mu)\ff{\pp}{\pp_{\xi(0)_i}},\ \     t\ge 0, \xi\in \C, \mu\in \scr P^\C.$$
Under reasonable conditions  on the coefficients $a_{ij}$ and $b_i$,
 the existence, uniqueness,   Lipschitz continuity in   Wasserstein distance, total variational norm and entropy, as well as  derivative estimates  are derived for the martingale solutions. \end{abstract} \noindent
 \vskip 2cm

\section{Introduction}

In this paper, we investigate nonlinear PDEs for probability measures on the  path space using path-distribution dependent SDEs. To explain the motivation of the study, let us start from the
 following classical   PDE on $\scr P(\R^d)$, the set of probability measures on $\R^d$ equipped with the weak topology:
\beq\label{E1} \pp_t \mu(t)= L^* \mu(t).\ \ t\ge 0,\end{equation} for a second-order differential operator
$$L:= \ff 1 2\sum_{i,j=1}^d a_{ij} \pp_i\pp_j+\sum_{i=1}^d b_i \pp_i,$$
where $a=(a_{ij}): \R^d\to \R^d\otimes \R^d$ and $b=(b_i): \R^d\to \R^d$ are locally integrable. \eqref{E1} is just the (linear) Fokker--Planck--Kolmogorov equation (FRKE) associated to the operator $L$ in the
sense of \cite{BKRS}. We call $\mu\in C(\mathbb R_+; \scr P(\R^d))$ a solution of \eqref{E1}, if
$$\int_{\R^d} f\d\mu(t)= \int_{\R^d} f\d\mu(0) +\int_0^t \d s\int_{\R^d}(Lf)\d\mu(s),\ \ t\ge 0, f\in C_0^\infty(\R^d).$$
To construct and analyze solutions of \eqref{E1} using the time marginal distributions of Markov processes as proposed by A. N. Kolmogorov \cite{KL}, K. It\^o  developed the theory of stochastic differential
equations (SDEs), see e.g.\cite{IW}.
Let $\si$ be a matrix-valued function such that $a=\si\si^*$, and let $W(t)$ be a $d$-dimensional Brownian motion. Consider the following It\^o SDE
\beq\label{E2} \d X(t)= b(X(t))\d t+ \si(X(t)) \d W(t).\end{equation}
By It\^o's formula, the time marginals $\mu(t):=\scr L_{X(t)}$ = the law of $X(t)$ for t $\geq 0,$ solve the equation \eqref{E1}. This enables one to investigate FPKEs using a probabilistic approach.

Obviously, \eqref{E1} is a linear equation. In applications, many important PDEs for probability measures (or probability densities) are nonlinear, see, for instance, \cite{CA, DV1, DV2, FG, Gu, V2} and
references within for the study of Landau type equations. Such PDEs are also of Fokker--Planck type, but are non-linear (see Sections 6.7 and 9.8 (v) in \cite{BKRS}).
To analyze non-linear FPKEs for probability measures, McKean-Vlasov equations are introduced by using SDEs with coefficients  depending on the distribution of the solution
  Consider the following distribution-dependent version of \eqref{E2}:
\beq\label{E3} \d X(t)= b(t,X(t),\scr L_{X(t)})\d t+ \si(t,X(t),\scr L_{X(t)}) \d W(t),\end{equation}
where $$b: \mathbb R_+\times \R^d\times\scr P(\R^d)\to \R^d,\ \ \si:  \mathbb R_+\times \R^d\times\scr P(\R^d)\to \R^d\otimes\R^d$$ are measurable.  For any $t\ge 0$ and $\mu\in \scr P(\R^d)$, consider
the second order differential operator
$$L_{t,\mu}:= \ff 1 2 \sum_{i,j=1}^d (\si\si^*)_{ij}(t,\cdot,\mu)\pp_i\pp_j +\sum_{i=1}^d b_i(t,\cdot,\mu)\pp_i.$$
Under reasonable integrability conditions on $\si$ and $b$, by It\^o's formula we see that for a solution $X(t)$ of \eqref{E3}, $\mu(t):=\scr L_{X(t)}$ solves the nonlinear FPKE
\beq\label{E4} \d \mu(t)= L_{t,\mu(t)}^* \mu(t)\end{equation} in the sense that
$$\int_{\R^d} f\d\mu(t)= \int_{\R^d} f\d\mu(0) +\int_0^t \d s\int_{\R^d}(L_{s,\mu(s)}f)\d\mu(s),\ \ t\ge 0, f\in C_0^\infty(\R^d).$$
   There are plentiful references on this type SDEs but most are concerning the existence, uniqueness and moments estimates, see  \cite{SZ, MV} and references within.
 In   the recent paper \cite{W16},  regularity estimates on the distribution, including  the exponential convergence  and gradient-Harnack type inequalities,    are presented. {\bf See also \cite{But,DV1, DV2, EGZ} for the study of ergodicity of distribution dependent SDEs.}

In the above two situations, the stochastic systems are Markovian (or memory-free); i.e. the evolution of the system does not depend on its past. However, many real-world models,
in particular those arising from  mathematical finance and biology, are with memory, so that  the associated evolution equations are path dependent, {\bf see, for instance, the monograph \cite{Moh} for specific models. See also  \cite{BWY17, BWY18, HMS} and references within for the study of regularities   and  ergodicity  for the associated FPKEs (i.e. the distribution of the functional solutions).}

In this paper, we investigate  nonlinear FPKEs on the path space by using path-distribution dependent SDEs.
In Section 2, we introduce the framework of the study and  the main  results on nonlinear FPKEs for probability measures  on  path space. To prove these results, we investigate the corresponding
path-distribution dependent SDEs in Sections 3-5, where strong/weak existence and uniqueness of solutions as well as Harnack type inequalities are derived respectively.  We will mainly follow the ideas
of \cite{W16}, but substantial additional efforts have to be made in order to generalize the results in there to the case, where the coefficients do not only depend on the time marginals,
but are also on the distribution of the path.

\section{Nonlinear PDEs for   measures on path space}

Throughout the paper, we fix $r_0>0$ and consider the path space $\C:= C([-r_0,0];\R^d)$ equipped with the uniform norm $\|\xi\|_\infty:=\sup_{\theta\in [-r_0.0]}|\xi(\theta)|. $ Let $\scr P_2^\C$
be the class of probability measures on $\C$ of finite second-order moment, i.e.  $\mu(\|\cdot\|_\infty^2):=\int_\C \|\xi\|_\infty^2\mu(\d\xi)<\infty.$ Then $\scr P_2^\C$ is a Polish space under
the Wasserstein distance
$$\W_2(\mu,\nu):=\inf_{\pi\in \C(\mu,\nu)} \bigg(\int_{\C\times\C} \|\xi-\eta\|_\infty^2\pi(\d\xi,\d\eta)\bigg)^{\ff 1 2},$$ where $\C(\mu,\nu)$ denotes the class of couplings for $\mu$ and $\nu$.
It is well known that  $(\scr P_2^\C, \W_2)$ is a Polish space and  the $\W_2$-metric  is consistent with the weak topology. We will study non-linear FPKEs on $\scr P_2^\C$.

Let \beq\label{*0D} b: \R_+\times \C\times \scr P_2^\C\to \R^d;\ \ \si: \R_+\times \C\times \scr P_2^\C\to \R^d\otimes\R^d\end{equation}
be measurable.
 For any $t\ge 0, \mu\in \scr P_2^\C$, consider the following differential operator $L_{t,\mu}$ from $C_0^\infty(\R^d)$ to the set of all $\B(\C)$-measurable functions: for $f\in C_0^\infty(\R^d)$,
 \beg{align*}
 (L_{t,\mu}f)(\xi):= \ff 1 2 \sum_{i,j=1}^d (\si\si^*)_{ij}(t, \xi, \mu) (\pp_i\pp_j f)(\xi(0))
 +\sum_{i=1}^d b_i(t,\xi,\mu) (\pp_if)(\xi(0)),\ \ \xi \in \C.
 \end{align*}
Then the associated  nonlinear FPKE for probability measures $(\mu_t)_{t\ge 0}$ on the path space $\C$ is
\beq\label{EM} \pp_t \mu(t) = L_{t,\mu_t}^* \mu_t,\end{equation}
where $\mu(t)$ is the marginal distribution of $\mu_t$ at $\theta=0$; i.e. $$\{\mu(t)\}(\d x):= \mu_t(\{\xi\in \C: \xi(0)\in \d x\}). $$
A continuous functional $\mu_\cdot: \R_+\to \scr P_2^\C$ is called a solution to \eqref{EM}, if $\int_0^t \d s\int_\C|L_{s,\mu_s}f|\d\mu_s<\infty$ for $f\in C_0^\infty(\R^d)$ and
\beq\label{S}\int_{\R^d} f\d\mu(t)=\int_{\R^d} f\d\mu(0) +\int_0^t \d s\int_\C(L_{s,\mu_s}f)\d\mu_s,\ \ t\ge 0,~ f\in C_0^\infty(\R^d).\end{equation}
{\bf Since \eqref{S} only characterizes the evolution of the marginal distribution $\mu(t)$ of $\mu_t$, the solution may be not unique. For instance, when $L_\mu$ is distribution independent and it has an invariant probability  measure $\mu(0)$, then any $\mu_t$ with marginal distribution $\mu(0)$ at $\theta=0$ solves \eqref{EM} in the sense of \eqref{S}. To select a   unique solution associated with the corresponding path-distribution dependent SDEs, 
we will only consider the martingale solutions of \eqref{EM}, which are a special class of solutions realized by  marginals of probability measures on the infinite-time path space $\C_\infty:=C([-r_0,\infty);\R^d)$. 
As shown in Theorem \ref{T1.1} below, in many cases the martingale solution is unique. }

For a probability
measure $ \mu^\infty$ on  $\C_\infty,$ consider  its marginal distributions
$$\mu^\infty(t):=  \mu^\infty \circ \{\pi(t)\}^{-1}\in \scr P(\R^d),\ \ \mu_t^\infty:=  \mu^\infty \circ\pi_t^{-1}\in \scr P(\C),\ \ t\ge 0,$$
where $\pi(t): C([-r_0,\infty);\R^d)\to\R^d$ and $\pi_t: C([-r_0,\infty);\R^d)\to\C$ are projection operators defined by
$$\pi(t) \xi= \xi(t)\in \R^d,\ \ \ \pi_t\xi= \xi_t\in \C \ \text{with}\ \xi_t(\theta):=\xi(t+\theta)\ \text{for}\  \theta\in [-r_0,0].$$

\beg{defn} A solution $(\mu_t)_{t\ge 0}$ of \eqref{EM} is called a martingale solution, if there exists a probability measure $ \mu^\infty$ on $\C_\infty$ such that
\beg{enumerate} \item[(1)]   $\mu_t=\mu_t^\infty$ for all $t\ge 0$.
\item[(2)] For any $f\in C_0^\infty(\R^d)$, the family of functionals
$$M^f(t):= f(\pi(t)\cdot)-\int_0^t (L_{s, \mu_s} f)(\pi_s\cdot)\d s,\ \ t\ge 0$$ on $\C_\infty$ is a $\mu^\infty$-martingale; that is,
$$\int_{A} M^f(t_2) \d\mu^\infty = \int_{A} M^f(t_1) \d\mu^\infty,\ \   t_2>t_1\ge 0,~ A\in \si (\pi(s): s\le t_1),$$ where $\si (\pi(s): s\le t_1)$ is  the $\si$-field on $\C_\infty$ induced by the projections $\pi(s)$ for $ s\in [-r_0,t_1]$.
\end{enumerate}  \end{defn}

To construct the martingale solutions of \eqref{EM} using path-distribution dependent SDEs, we need the following assumptions.

\beg{enumerate} \item[$(H1)$] (Continuity)  For every $t\ge 0$, $b(t,\cdot,\cdot)$ is continuous on $\C \times\scr P_2^\C$, and  there exist locally bounded functions $\aa_1,\aa_2: \R_+\to \R_+$  such that
 $$\|\si(t,\xi,\mu)- \si(t,\eta,\nu)\|^2\le \aa_1(t)\|\xi-\eta\|_{\infty}^2+ \aa_2(t) \W_2(\mu,\nu)^2,\ \ t\ge 0; \xi,\eta\in \C; \mu,\nu\in
\scr P_2^\C.$$
 \item[$(H2)$] (Monotonicity) There exist a constant $\kk\ge 0$ and  locally bounded functions   $\bb_1,\bb_2: \R_+\to \R_+$ such that
 \beg{align*} &2\langle  b(t,\xi,\mu)- b(t,\eta,\nu), \xi(0)-\eta(0)\rangle+\|\si(t,\xi,\mu)- \si(t,\eta,\nu)\|_{HS}^2\\
 &\le \bb_1(t) \|\xi-\eta\|_{\infty}^2+ \bb_2(t)\W_2(\mu,\nu)^2 - \kk |\xi(0)-\eta(0)|^2,\ \ t\ge 0; \xi,\eta\in \C; \mu,\nu\in
\scr P_2^\C.\end{align*}
\item[$(H3)$] (Growth)  $b$ is bounded on bounded sets in $[0,\infty)\times \C\times \scr P_2^\C$, and there exists a locally bounded function  $K: \R_+\to\R_+   $ such that
$$|b(t,0,\mu)|^2+ \|\si(t,0,\mu)\|^{2}\le K (t) \big\{1+\mu(\|\cdot\|_{\infty}^2)\big\},\ \ t\ge 0,~ \mu\in \scr P_2^\C.$$
\end{enumerate}

The following result characterizes the martingale solutions of \eqref{EM} with $\W_2$-Lipschitz estimate.

\beg{thm}\label{T1.1} Assume $(H1)$-$(H3)$.  Then for any $\mu_0\in \scr P_2^\C$,  there exists a unique martingale solution   $ (\mu_t)_{t\ge 0} $ of $\eqref{EM}$. Moreover,
 \beg{enumerate} \item[$(1)$] $\mu_t( \|\cdot\|_\infty^2)$ is locally bounded in $t$.
\item[$(2)$] For any  two martingale solutions $(\mu_t)_{t\ge 0}$ and $ (\nu_t)_{t\ge 0}$ of $\eqref{EM}$,
\begin{align*} &\W_2(\mu_t,  \nu_t)^2\le \inf_{\vv \in [0,1]} \left\{\ff{ \W_2(\mu_0,\nu_0)^2   }{1-\vv} \right.\\
   &\left. \times \inf_{\dd\in [0,\kk]} \exp\bigg[(r_0-t)\dd+ \ff{\e^{\dd r_0}}{1-\vv}\int_0^t\Big\{\ff{4  (\aa_1(r)+\aa_2(r))}\vv +  \bb_1(r)+\bb_2(r) \Big\}\d r\bigg] \right\} \end{align*}
holds for all $t\ge 0$ and $\vv\in (0,1)$.
\end{enumerate} \end{thm}

From now on, for any $\nu_0,\mu_0\in \scr P_2^\C$, we denote $\mu_t$ and $\nu_t$ the martingale solutions of \eqref{EM} staring at $\mu_0$ and $\nu_0$ respectively.

 To estimate the continuity of $\mu_t$ in $\mu_0$ with respect to entropy and total variational norm, we  make the following stronger assumption.
 \beg{enumerate} \item[{\bf(A)}] $\si(t,x)$ is invertible,  and  there exist  increasing functions  $\kk_0,\kk_1,\kk_2, \ll: \R_+\to  \R_+$ such that for any $t\ge 0, x,y\in\R^d,\xi,\eta\in \C$
 and $\mu,\nu\in \scr P_2^\C$,
 \beg{align*}
& |b(t,0,\mu)|^2+\|\si(t,x)\|^2\le \kk_0(t)(1+|x|^2+\mu(\|\cdot\|_{\infty}^2)),\\
  & \|\si(t,\cdot)^{-1}\|_\infty\le \ll(t),\ \ \|\si(t,x)-\si(t,y)\|_{HS}^2 \le \kk_1(t)|x-y|^2,\\
 &|b(t,\xi,\mu)-b(t,\eta,\nu)|\le\kk_{2}(t)(\|\xi-\eta\|_{\infty}+\W_2(\mu,\nu)).\end{align*}
\end{enumerate}
   Recall that for any two probability measures $\mu,\nu$ on   some measurable space $(E,\scr F)$, the entropy and variational norm are defined as follows:
$$\Ent(\nu|\mu):= \beg{cases} \int (\log \ff{\d\nu}{\d\mu})\d\nu, \ &\text{if}\ \nu\ \text{ is\ absolutely\ continuous\ with\ respect\ to}\ \mu,\\
 \infty,\ &\text{otherwise;}\end{cases}$$ and
$$\|\mu-\nu\|_{var} := \sup_{A\in\F}|\mu(A)-\nu(A)|.$$ By Pinsker's inequality (see \cite{CK, Pin}),
\beq\label{ETX} \|\mu-\nu\|_{var}^2\le \ff 1 2 \Ent(\nu|\mu),\ \ \mu,\nu\in \scr P(E).\end{equation}   Then \eqref{TT2} below implies
\beq\label{TT}\|\mu_t-\nu_t\|_{var}^2\le  \ff{\psi(t)}{2(t-r_0)}\W_2(\mu_0,\nu_0)^2 ,\ \ t>r_0,\ \mu_0,\nu_0\in \scr P_2^\C,  \end{equation}
for some $\psi\in C(\R_+;\R_+)$.
  There are a lot of examples where $\W_2(\mu_n,\mu_0)\to 0$ but $\mu_n$ is singular with respect to $\mu_0$ such that $\Ent(\mu_n|\mu_0)=\infty$ and $\|\mu_n-\mu_0\|_{var}=1.$ So,
  both \eqref{TT} and \eqref{TT2}  are non-trivial. Indeed, these estimates correspond to  the log-Harnack inequality  for the associated semigroups, see Theorem \ref{T3.1} below for details.

\beg{thm}\label{T1.2} Assume {\bf (A)}.  \beg{enumerate} \item[$(1)$] There exists   $\psi\in C(\R_+;\R_+)$ such that
\beq\label{TT2} \Ent(\nu_t|\mu_t)  \le \ff{\psi(t)}{t-r_0}\W_2(\mu_0,\nu_0)^2 ,\ \ t>r_0, \mu_0,\nu_0\in \scr P_2^\C. \end{equation}
 \item[$(2)$] If there exists an increasing   function $\kk_3: \R_+\to \R_+$ such that
 \beq\label{*P} \|\si(t,x)-\si(t,y)\| \le \kk_3(t)(1\land |x-y|), \ \ t\ge 0, x,y\in\R^d,\end{equation}  then there exists a positive continuous function $H$ defined on the domain
 $$D:=\{(p,t): t\ge 0, p>(1+\kk_3(t)\ll(t))^2\},$$ such that
 $$\int_\C \Big(\ff{\d\nu_t}{\d\mu_t}\Big)^{\ff 1 p} \d\nu_t\le \inf_{\pi\in\scr C(\mu_0,\nu_0)} \int_{\C\times\C}  \e^{H(p,t)\big(1+\ff{|\xi(0)-\eta(0)|^2}{t-r_0} + \|\xi-\eta\|_\infty^2\big)}\d\pi$$
 holds for all
 $t>r_0$ and $p>(1+\kk_3(t)\ll(t))^2.$
   \end{enumerate} \end{thm}

\paragraph{Remark 2.1.}  According to Theorem \ref{T1.1}(2), if there exists a constant $\vv\in (0,1)$ such that
\beq\label{EXO} \limsup_{t\to\infty}  \ff 1 t \int_0^t \Big(\ff{4(\aa_1(s)+\aa_2(s))}{\vv(1-\vv)} + \ff{\bb_1(s)+\bb_2(s)}{1-\vv}\Big)\d s <\sup_{\dd\in [0,\kk]}\dd\e^{-\dd r_0},\end{equation}
then \beq\label{EXO2} \W_2(\mu_t,  \nu_t)^2\le   c\e^{-\ll t} \W_2(\mu_0,\nu_0)^2,\ \ t\ge 0, \end{equation} holds for some constants $c,\ll>0$;
i.e.\,the solution to \eqref{EM} has  exponential contraction
in $\W_2$. If
$\si(t,\cdot,\cdot)$ and $b(t,\cdot,\cdot)$ do not depend on $t$, i.e.\,the equation is time-homogenous, we  $\mu_t= P_t^*\mu_0$. By the uniqueness we see that $P_t^*$ is a semigroup, i.e.\,$P^*_{t+s}=P_t^*P_s^*, s,t\ge 0.$ Then \eqref{EXO}  implies that   $P_t^*$  has a unique  invariant probability measure $\mu\in \scr P_2^\C$.
Combining  \eqref{EXO2} with the semigroup property  of $P_t^*$ and \eqref{TT}-\eqref{TT2}, we conclude that \eqref{EXO} also implies the exponential convergence in entropy and total variational norm:
$$\max\{\Ent(\nu_t|\mu), \|\mu-\nu_t\|^2_{var}\} \le c_1 \W_2(\mu,\nu_{t-1})^2\le c_2\e^{-\ll t} \W_2(\mu,\nu_0)^2,\ \ t\ge 1, \nu_0\in\scr P_2^\C$$ for some constants $c_1,c_2>0$.

\

   Finally, we investigate the shift quasi-invariance and differentiability of $\mu_t$ along Cameron--Martin vectors in $\H^1:= \{\xi\in \C: \int_{-r_0}^0 |\xi'(s)|^2\d s<\infty\}.$  For  $\eta\in \C$ and a
   probability measure $\mu$ on $ \C$, we say that $\mu$ is differentiable along $\xi$ if for any $A\in \B(\C)$, $\pp_\xi \mu(A):= \ff{\d}{\d \vv} \mu(A+\vv \xi)\big|_{\vv=0}$ exists and $\pp_\xi\mu(\cdot)$ is a
   signed measure on $\C$.

  \beg{thm}\label{T1.3} Assume {\bf (A)} and let $b(t,\cdot,\mu)$ be differentiable on $\C$, $\si(t,x)=\si(t)$ be independent of $x$. Then for any $t>r_0, \eta\in \H^1$ and $\mu_0\in \scr P_2^\C$,
  $\mu_t$ is differentiable along $\eta$, both $\pp_\eta\mu_t$ and $\mu_t(\cdot+\eta)$ are absolutely continuous with respect to $\mu_t$, and for some $\Psi\in C(\R_+;\R_+)$
  \beg{align*} &\int_\C \Big(\log \ff{\d\mu_t(\cdot+\eta)}{\d\mu_t}\Big)\d\mu_t (\cdot+\eta) \le \Psi(t) \Big(\ff{|\eta(-r_0)|^2}{t-r_0} + \|\eta\|_{\H^1}^2\Big),\\
  &\int_\C \Big( \ff{\d\mu_t(\cdot+\eta)}{\d\mu_t}\Big)^{\ff 1 p} \d\mu_t (\cdot+\eta)\le \exp\bigg[\Psi(t) \Big(\ff{|\eta(-r_0)|^2}{t-r_0} + \|\eta\|_{\H^1}^2\Big)\bigg],\ \ p>1,\\
  & \int_\C\Big| \ff{\d\pp_\eta\mu_t}{\d\mu_t}\Big|^2 \d\mu_t \le \Psi(t) \Big(\ff{|\eta(-r_0)|^2}{t-r_0} + \|\eta\|_{\H^1}^2\Big).\end{align*}
   \end{thm}

\beg{proof} [Proof of Theorems  \ref{T1.1}-\ref{T1.3}]
For $\mu_0\in \scr P_2^\C$, take a $\F_0$-measurable random variable $X_0$ on $\C$ such that $\scr L_{X_0}=\mu_0$. According to Theorem \ref{T2.1}, Corollary \ref{C3.3}, Corollary \ref{C4.2} and \eqref{ETX},
$\mu_t:= \scr L_{X_t}$ satisfies the estimates in Theorems \ref{T1.1}-\ref{T1.3} under the corresponding assumptions. So, it suffices to show that  $(\scr L_{X_t})_{t\ge 0}$ is
  the unique martingale solution of \eqref{EM}.

Let  $\mu^\infty = \L_{\{X(s)\}_{s\in [-r_0.\infty)}}$. We have $ \scr L_{X_t} =\mu_t^\infty$. By \eqref{ED1} and It\^o's formula, for any $f\in C_0^\infty(\R^d)$,  $(M^f(t))_{t\ge 0}$ is a $\mu^\infty$-martingale
such that $\mu_t:=\scr L_{X_t}$ satisfies
\beg{align*} \int_{\R^d} f\d\mu(t)&=\E f(X(t))= \E f(X(0) ) +\int_0^t \E (L_{s,\mu_s}f)(X_s)\d s\\
&= \int_{\R^d} f\d\mu(0) + \int_0^t\d s\int_{\C} (L_{s,\mu_s}f)\d\mu_s,\ \ t\ge 0,  f\in C_0^\infty(\R^d).\end{align*}
Therefore, $\L_{\{X(s)\}_{s\in [-r_0.\infty)}}$ is a martingale solution  of \eqref{EM}.
When the coefficients are distribution-free, it is well known that the weak solution of \eqref{ED1} is equivalent to the martingale solution, so that the uniqueness of the martingale solutions of \eqref{EM}
follows   from Theorem \ref{T2.1}(3) below. In the following,  we explain that the same is true for the present distribution dependent case.

Let $\mu_t=\mu_t^\infty$, for some  probability measure $\mu^\infty$ on $\C_\infty$, be a martingale solution of \eqref{EM}.  We intend to prove $\mu^\infty= \L_{\{X(s)\}_{s\in [-r_0.\infty)}},$ so that
the martingale solution is unique.
Let   $\bar \OO:=\C_\infty, \bar \F_t$ for $t\ge 0$ be the completion of $\si(\pi(s): s\le t)$ with respect to $\mu^\infty$, and $\bar \P:=\mu^\infty$.
By Theorem \ref{T2.1}(3) below, it suffices to prove that the
coordinate process
$$\bar X (t)(\oo):= \oo(t),\  \ t\ge 0, \oo\in \bar\OO$$ is a weak solution to \eqref{ED1}.
To this end, for the given $(\mu_t)_{t\ge 0}$,  define
$$\bar\si(t,\xi):= \si(t,\xi,\mu_t),\ \ \ \bar b(t,\xi)= b(t,\xi,\mu_t),\ \ t\ge 0,~ \xi\in\C,$$ and consider the corresponding operator
$$(\bar L_t f)(\xi):= \ff 1 2\sum_{i,j=1}^d (\bar\si\bar\si^*)_{ij}(t,\xi)(\pp_i\pp_jf)(\xi(0)) +\sum_{i=1}^d \bar b_i(t,\xi) (\pp_if)(\xi(0)),\ \ t\ge 0, \xi\in \C$$ for $f\in C_0^\infty(\R^d).$
Since $(\mu_t)_{t\ge 0}$ is a martingale solution of \eqref{EM}, for any $f\in C_0^\infty(\R^d)$, the process
$$  M^f(t):=f(\bar X(t))-f(\bar X(0))-\int_0^t (\bar L_sf)(\bar X_s)\d s,\ \ t\ge 0$$ is a martingale on the probability space $(\bar\OO, (\bar\F_t)_{t\ge 0}, \bar\P).$
 By $(H1)$-$(H3)$,  the martingale property also holds for $f$ being polynomials of order 2.
 In particular,    by taking $f(x)= x$ we see that
\beq\label{AB} M(t):=M^f(t)= \bar X(t)-\bar X(0)-\int_0^t \bar b(s, \bar X_s)\d s\end{equation} is a $\R^d$-valued  martingale, and with $f(x):=x_ix_j$ we conclude that
$$\<M_i,M_j\>(t) =\int_0^t (\bar\si\bar\si^*)_{ij}(s, \bar X_s)\d s,\ \ 1\le i,j\le d.$$
Then according to Stroock--Varadhan   (see, for example,  Theorems 4.5.1 and 4.5.2 in \cite{SV}),
we may construct   a  $d$-dimensional Brownian motion $\tt W(t)$ on a product probability space of $(\tt \OO, \tt\F_t,\tt \P)$ with $(\bar \OO, \bar\F_t,\bar \P)$ as a marginal space, and when $\si$ is
invertible these two spaces coincide,  such that
$$M(t)=\int_0^t \bar\si(s,\bar X_s)\d \tt W(s),\ \ t\ge 0.$$
Combining this with \eqref{AB}, we see that $\bar X(t)$ solves the
 stochastic functional differential equation
\beq\label{DD0} \d \bar X(t)= \bar b(t,\bar X_t)\d t+\bar \si(t,\bar X_t)\d \tt W(t)\end{equation}  with $\L_{\bar X_0}|_{\tt\P}=\L_{\bar X_0}|_{\bar\P}=\mu_0.$
Since, by definition, $\mu_t=\scr L_{\bar X_t}|_{\bar \P}=\scr L_{\bar X_t}|_{\tt \P}$,
$\bar X(t)$ solves   the path-distribution dependent SDE
$$\d \bar X(t)= b(t,\bar X_t, \L_{\bar X_t}|_{\tt\P}) \d t+ \si(t,\bar X_t, \L_{\bar X_t}|_{\tt\P})\d \tt W(t),$$ i.e.\,$(\bar X, \tt W)$ is a weak solution of \eqref{ED1}.
Noting that $\mu^\infty:= \scr L_{\bar X|\bar\P} =\scr L_{\bar X|\tt\P},$  by the weak uniqueness of \eqref{ED1} due to Theorem \ref{T2.1}(3) below,  we obtain $\mu^\infty = \L_{\{X(s)\}_{s\in [-r_0,\infty)}}$
as desired.
\end{proof}

\section{Path-distribution dependent SDEs}

Recall that for $\gg(\cdot)\in C([-r_0,\infty);\R^d)$, the segment functional $\gg_\cdot\in C(\R_+;\C)$ is defined by
$$\gg_t(\theta):= \gg(t+\theta),\ \ \theta\in [-r_0,0], t\ge 0.$$ For $\si,b$ in \eqref{*0D},
consider the following path-distribution dependent SDE on $\R^d$:
\beq\label{ED1}
\d X(t)= b(t,X_t,\L_{X_t})\,\d t+ \sigma(t,X_t, \L_{X_t})\,\d W(t),
\end{equation}
where $W=(W(t))_{t\geq 0}$ is a $d$-dimensional standard Brownian motion with respect to a complete filtered probability space $(\OO, \F, \{\F_{t}\}_{t\ge 0}, \P)$, $\L_{X_t}$ is the distribution of $X_t$.
We investigate the strong solutions of \eqref{ED1} and determine properties, of their distributions.

We first recall the definition of the strong and weak solutions, see for instance \cite[Definition 1.1]{W16} in the path independent setting.
For simplicity, we will only consider square integrable solutions.
\beg{defn}  $(1)$ For any $s\ge 0$, a continuous adapted process $(X_{s,t})_{t\ge s}$ on $\C$ is called a (strong) solution of \eqref{ED1} from time $s$, if
  $$\E\|X_{s,t}\|_\infty^2 +\int_s^t \E\big\{|b(r,X_{s,r},\L_{X_{s,r}})|+\|\si(r,X_{s,r}, \L_{X_{s,r}})\|^2\big\}\d r<\infty,\ \ t\ge s,$$   and  $(X_s,(t):= X_{s,t}(0))_{t\ge s}$ satisfies $\P$-a.s.
$$X_s,(t) = X_s(s) +\int_s^t b(r,X_{s,r}, \L_{X_{s,r}})\d r + \int_s^t \si(r,X_{s,r},\L_{X_{s,r}})\d W(r),\ \ t\ge s.$$
We say that \eqref{ED1} has (strong or pathwise) existence and uniqueness, if for any $s\ge 0$ and $\F_s$-measurable random variable $X_{s,s}$ with $\E\|X_{s,s}\|_{\infty}^2<\infty$, the equation from
time $s$ has a unique solution  $(X_{s,t})_{t\ge s}$. When $s=0$ we    simply denote $X_{0,}=X$; i.e.\,$X_{0,}(t)=X(t), X_{0,t}=X_t, t\ge 0$.

$(2)$ A couple $(\tt X_{s,t}, \tt W(t))_{t\ge s}$ is called a weak solution to \eqref{ED1} from time $s$, if $\tt W(t)$  is a $d$-dimensional Brownian motion a complete filtered probability space
$ (\tt\OO, \{\tt\F_t\}_{t\ge s}, \tt\P)$, and   $\tt X_{s,t}$ solves
\beq\label{E1'} \d \tt X_{s,}(t)= b(t,\tt X_{s,t},  \L_{\tt X_{s,t}}|_{\tt\P})\d t + \si(t,\tt X_{s,t},  \L_{\tt X_{s,t}}|_{\tt\P})\d \tt W(t),\ \ t\ge s.\end{equation}

$(3)$  \eqref{ED1} is said to satisfy weak uniqueness, if for any $s\ge 0$, the distribution of a weak solution $(X_{s,t})_{t\ge s}$ to \eqref{ED1} from $s\ge 0$ is uniquely determined by $\scr L_{X_{s,s}}$.
\end{defn}

 When \eqref{ED1} has strong existence and uniqueness, the solution $(X_t)_{t\ge 0}$ is a Markov process in the sense that for any $s\ge 0$,  $(X_t)_{t\ge s}$ is determined by solving the equation
 from time $s$ with initial state $X_s$. More precisely, letting $\{X_{s,t}^{\xi}\}_{t\ge s}$ denote the solution of the equation from time $s$ with initial state  $X_{s,s}=\xi$,
 the existence and uniqueness imply
 \beq\label{MK}
 X_{s,t}^{\xi}= X_{u,t}^{X_{s,u}^{\xi}},\ \ t\ge u\ge  s\ge 0, \xi\ {\rm is}\ \F_s\text{-measurable\ with\ } \E\|\xi\|_{\infty}^2<\infty.
 \end{equation}

When \eqref{ED1}  also has   weak uniqueness, we may define a semigroup $(P_{s,t}^*)_{t\ge s}$ on $\scr P_2^\C$ by letting $P_{s,t}^*\mu=\L_{X_{s,t}}$ for   $\L_{X_{s,s}}=\mu\in \scr P_2^\C$.
Indeed, by \eqref{MK} we have
\beq\label{SM}
P_{s,t}^*= P_{u,t}^* P_{s,u}^*,\ \ t\ge u\ge  s\ge 0.
\end{equation}
For simplicity we set $P_t^\ast=P_{0,t}^\ast,~ t\ge 0$.

 \beg{thm}\label{T2.1} Assume $(H1)$-$(H3)$.
\beg{enumerate} \item[$(1)$] For any $s\geq 0$ and $X_{s,s}\in L^2(\OO\to \C;\F_s)$, $\eqref{ED1}$ has a unique strong solution $(X_{s,t})_{t\ge s}$ with
\beq\label{ES1}
\E \sup_{t\in [s,T]}\|X_{s,t}\|_{\infty}^{2}\le H(T) (1+ \E\|X_{s,s}\|_{\infty}^{2}),\ \ T\ge t\ge s\ge  0 \end{equation}   for some increasing function $H:\R_+\to \R_+$.
\item[$(2)$] For any two solutions $X_{s,t}$ and $Y_{s,t}$ of $\eqref{ED1}$ with $\L_{X_{s,s}},\L_{Y_{s,s}}\in \scr P_2^\C$,
\beg{align*}
&\E \|X_{s,t}-Y_{s,t}\|_{\infty}^2\le  \inf_{\vv\in (0,1)} \left\{ \ff{ \E\|X_{s,s}-Y_{s,s}\|_{\infty}^2   }{1-\vv} \right. \\
&\left.\times \inf_{\dd\in [0,\kk],\vv\in (0,1)} \exp\bigg[(r_0+s-t)\dd + \ff{\e^{\dd r_0}}{1-\vv}\int_s^t\Big\{\ff{4  (\aa_1(r)+\aa_2(r))}\vv +  \bb_1(r)+\bb_2(r) \Big\}\d r\bigg] \right\}.\end{align*}
\item[$(3)$]   $\eqref{ED1}$ satisfies weak uniqueness,  and for any $t\ge 0$,
\beg{align*} &\W_2(P_t^*\mu_0,  P_t^*\nu_0)^2\le \inf_{\vv \in (0,1)} \left\{ \ff{ \W_2(\mu_0,\nu_0)^2   }{1-\vv} \right. \\
&\times \inf_{\dd\in [0,\kk], \vv\in (0,1)}\exp\bigg[(r_0-t)\dd+\ff{\e^{\dd r_0}}{1-\vv}\int_0^t\Big\{\ff{4  (\aa_1(r)+\aa_2(r))}\vv +  \bb_1(r)+\bb_2(r) \Big\}\d r\bigg]. \end{align*}
\end{enumerate} \end{thm}

We will prove this result by using the argument of \cite{W16}. For fixed  $s\ge 0$ and $\F_s$-measurable  $\C$-valued random variable $X_{s,s}$ with $\E\|X_{s,s}\|_{\infty}^2<\infty$, we construct the solution of \eqref{ED1} by iterating in distribution as follows. Firstly, let
$$X^{(0)}_{s,t}(\theta)=X_{s,s}\big(0\land(t-s+\theta)\big)\ \text{for}\ \theta\in [-r_0,0],     \ \mu_{s,t}^{(0)}=\L_{X^{(0)}_{s,t}},\ \ t\ge s.$$  For any $n\ge 1,$ let $(X_{s,t}^{(n)})_{t\ge s}$ solve the classical path-dependent SDE
\beq\label{EN}
\d X^{(n)}_{s,}(t)= b(t,X^{(n)}_{s,t}, \mu_{s,t}^{(n-1)}) \d t + \si(t,X^{(n)}_{s,t},\mu_{s,t}^{(n-1)})\,\d W(t),\ \ X_{s,s}^{(n)}=X_{s,s}, t\ge s,
\end{equation}
where $\mu_{s,t}^{(n-1)}:=\L_{X_{s,t}^{(n-1)}}$  and $X^{(n)}_{s,t}(\theta):= X^{(n)}_{s,}(t-s+\theta)$ for $\theta\in [-r_0, 0]$.

\beg{lem} \label{L2.1} Assume $(H1)$-$(H3)$. For every $n\ge 1$, the path-dependent SDE $\eqref{EN}$ has a unique strong solution $X^{(n)}_{s,t}$ with
\beq\label{*2} \E\sup_{t\in [s-r_0,T]} |X^{(n)}_{s,}(t)|^2<\infty,\ \ T>s, n\ge 1.\end{equation} Moreover, for any $T>0$, there exists $t_0>0$ such that for all $s\in [0,T]$ and $X_{s,s}\in L^2(\OO\to\C;\scr F_s)$,
\beq\label{*3}
\E \sup_{ t\in [s, s+t_0]} |X^{(n+1)}_{s,}(t)-X^{(n)}_{s,}(t)|^2\le 4 \e^{-n}  \E\sup_{t\in [s,s+t_0]} |X^{(1)}_{s,}(t)|^2,\ \ s\in [0,T], n\ge 1.
\end{equation}
\end{lem}

\beg{proof} The proof is similar to that of \cite[Lemma 2.1]{W16}.
Without loss of generality, we may assume that  $s=0$ and simply denote  $X_{0,}(t)=X(t), X_{0,t}=X_t, t\ge 0$.

(1) We first prove that  the SDE \eqref{EN} has a unique strong solution and \eqref{*2} holds.

For   $n=1$, let
$$\bar b(t,\xi)= b(t,\xi, \mu_{t}^{(0)}),\ \ \bar\si(t,\xi)= \si (t,\xi, \mu_{t}^{(0)}),\ \ t\ge 0, \xi\in\C.$$ Then \eqref{EN} reduces to
\beq\label{EN*} \d X^{(1)}(t)= \bar b(t, X_{t}^{(1)})\d t + \bar\si(t, X_{t}^{(1)})\d W(t),\ \ X_{0}^{(1)}=X_{0}, t\ge 0.\end{equation}
By $(H1)$-$(H3)$, the coefficients $\bar b$ and $\bar \si$ satisfy the standard monotonicity condition which imply strong existence, uniqueness and non-explosion for
the stochastic functional differential equation \eqref{EN*}, see e.g. \cite[Corollary 4.1.2]{Wbook} with $D=\mathbb{R}^d$ and $u_n=1$. It is also standard to prove  \eqref{*2} using It\^o's formula
\beg{align*} \d |X^{(1)}(t)|^2  &=2\big\<\si(t,X_t^{(1)},\mu_t^{(0)})\d W(t), X^{(1)}(t)\big\>\\
&\quad + \big\{ 2\big\<b(t, X_t^{(1)},  \mu_t^{(0)}), X^{(1)}(t)\big\> + \|\si(t,X_t^{(1)},\mu_t^{(0)})\|_{HS}^2\big\}\d t.\end{align*}
By  $(H1)$-$(H3)$,     there exists an increasing function $H:\R_+\to\R_+$ such that
\beg{align*}& 2\big\<b(t,\xi,\mu_t^{(0)}), \xi(0)\big\> + \|\si(t,\xi,\mu_t^{(0)})\|_{HS}^2\\
&\le 2\big\<b(t,\xi,\mu_t^{(0)})-b(t, 0,\mu_t^{(0)}), \xi(0) \big\>
  +   2|b(t,0,\mu_t^{(0)})|\cdot |\xi(0)|\\
  &+   2\|\si(t,\xi,\mu_t^{(0)})-\si(t,0,\mu_t^{(0)})\|^2_{HS}+2\|\si(t,0,\mu_t^{(0)})\|^2_{HS}\\
&\le H(t) \big\{1+ \|\xi\|_{\infty}^2+\mu_t^{(0)}(\|\cdot\|_{\infty}^2)\big\},\ \ \ t\ge 0, \xi\in \C.\end{align*}
Combining this with $(H3)$ and applying the BDG inequality for $p = 1$, for any $N\in [1,\infty)$ and $\tau_N:= \inf\{t\ge 0: |X^{(1)}(t)|\ge N\}$, we have
\beg{align*} &\E \sup_{s\in [-r_0,t\land \tau_N]} |X^{(1)}(s)|^2
 \le 4\mathbb{E}\|X^{(1)}_0\|_{\infty}^2+2H(t) \E\int_0^{t\land \tau_N} \big(1+\|X_s^{(1)}\|_{\infty}^2 + \mu_s^{(0)}(\|\cdot\|_{\infty}^2)\big)\d s \\
&\qquad + 4 H(t) \E\bigg(\int_0^{t\land \tau_N} |X^{(1)}(s)|^2 \big(1+ \|X_s^{(1)}\|_{\infty}^2 + \mu_s^{(0)}(\|\cdot\|_{\infty}^2)\big)\d s\bigg)^{\ff 1 2}\\
&\le 4\mathbb{E}\|X^{(1)}_0\|_{\infty}^2+ \ff 1 2 \E \sup_{s\in [-r_0,t\land \tau_N]} |X^{(1)}(s)|^2\\
&+\{2H(t)+ 8H(t)^2\} \E  \int_0^{t\land \tau_N} \big(1+\|X_s^{(1)}\|_{\infty}^2 + \mu_s^{(0)}(\|\cdot\|_{\infty}^2) \big)\d s,\ \ t\ge 0.\end{align*}
This implies
\begin{equation*}\begin{split}
&\E \sup_{s\in [-r_0,t\land \tau_N]} |X_s^{(1)}|^2 \le 8\mathbb{E}\|X^{(1)}_0\|_{\infty}^2\\
&+ \{4H(t)+ 16H(t)^2\} \int_0^{t} \big\{1+ \E \sup_{r\in [-r_0,s\land \tau_N]}|X^{(1)}(r)|^2 + \mu_s^{(0)}(\|\cdot\|_{\infty}^2)\big\}\d s,\ \ t\ge 0. \end{split}\end{equation*}
By first applying Gronwall's Lemma then letting $N\to\infty$, we arrive at
$$\E \sup_{s\in [-r_0,t]} |X^{(1)}(s)|^2<\infty,\ \ t\ge 0.$$
Therefore,   \eqref{*2}  holds  for $n=1$.

Now, assuming that the assertion holds for $n=k$ for some $k\ge 1$, we intend to prove it for $n=k+1$. This can be done by repeating the above argument with
 $(X_\cdot^{(k+1)}, \mu_\cdot^{(k)}, X_\cdot^{(k)})$ replacing $(X_\cdot^{(1)}, \mu_\cdot^{(0)},X_\cdot^{(0)})$, so, we omit the proof.

(2) To prove \eqref{*3}, let
\beg{align*} &\xi^{(n)}(t) = X^{(n+1)}(t)- X^{(n)}(t),\\
&\LL_t^{(n)}= \si(t,X_t^{(n+1)},\mu_t^{(n)})- \si(t,X_t^{(n)},\mu_t^{(n-1)}),\\
&B_t^{(n)}= b(t,X_t^{(n+1)}, \mu_t^{(n)} ) - b(t,X_t^{(n)}, \mu_t^{(n-1)} ).\end{align*}
 By $(H2)$ and It\^o's formula,  there exists an increasing function $K_1:\R_+\to\R_+$ such that
$$\d |\xi^{(n)}(t)|^2  \le 2 \<\LL_t^{(n)} \d W(t), \xi^{(n)}(t)\>+  K_1(t) \big\{\|\xi_t^{(n)}\|_{\infty}^2 + \W_2(\mu_t^{(n)}, \mu_t^{(n-1)})^2\big\}\d t.$$
By the BDG inequality for $p = 1$ and since $\W_2(\mu_s^{(n)}, \mu_s^{(n-1)})^2\le \E\|\xi_s^{(n)}\|_\infty^2$, we obtain
 \beg{align*}& \E \sup_{s\in [0,t]}  |\xi^{(n)}(s)|^2\le 2 \E\sup_{s\in [0,t]} \int_0^s \<\LL_r^{(n)}\d W(r), \xi^{(n)}(r)\>\\
  &\qquad +  K_1(t) \int_0^t \Big\{\E   \|\xi_s^{(n)}\|_{\infty}^2 + \W_2(\mu_s^{(n)}, \mu_s^{(n-1)})^2\Big\} \d s\\
  & \le   4\E\bigg(\int_0^t \big\{|\xi^{(n)}(s)|^2 \|\LL_s^{(n)}\|^2\big\}\d s\bigg)^{\ff 1 2}
 + K_1(t) \int_0^t \Big\{\E   \|\xi_s^{(n)}\|_{\infty}^2 + \W_2(\mu_s^{(n)}, \mu_s^{(n-1)})^2\Big\} \d s\\
 &\le \ff 1 2\E \sup_{s\in [0,t]}  |\xi^{(n)}(s)|^2 + 8\int_0^t\E \|\LL_s^{(n)}\|^2\d s + K_1(t) \int_0^t \Big\{\E   \|\xi_s^{(n)}\|_{\infty}^2 + \W_2(\mu_s^{(n)}, \mu_s^{(n-1)})^2\Big\} \d s. \end{align*}
 Combining this and $(H1)$ we deduce that
$$\E \sup_{s\in [0,t]}  |\xi^{(n)}(s)|^2 \le  K_2(t) \int_0^t \Big\{\E   \sup_{r\in[0,s]}|\xi^{(n)}(r)|^2 + \W_2(\mu_s^{(n)}, \mu_s^{(n-1)})^2\Big\} \d s,\ \ t\ge 0 $$ for some increasing function
$K_2: \R_+\to \R_+.$  By
  Gronwall's Lemma, we obtain
\beg{align*} & \E \sup_{s\in [0,t]}  |\xi^{(n)}(s)|^2 \le tK_2(t) \e^{tK_2(t)} \sup_{s\in [0,t]} \W_2(\mu_s^{(n)}, \mu_s^{(n-1)})^2\\
&\le tK_2(t) \e^{tK_2(t)}
  \E  \sup_{s\in [0,t]}|\xi^{(n-1)}(s)|^2,\ \  \ t\ge 0.\end{align*}
 Taking $t_0>0$ such that $  t_0K_2(T) \e^{t_0K_2(T)}\le \e^{-1}$, we arrive  at
 $$\E \sup_{s\in [0,t_0]} |\xi^{(n)}(s)|^2\le \e^{-1}\E \sup_{s\in [0,t_0]} |\xi^{(n-1)}(s)|^2,\ \ n\ge 1.$$ Since
 $$\E \sup_{s\in [0,t_0]} |\xi^{(0)}(s)|^2\le 2  \E  \Big\{|X(0)|^2+ \sup_{s\in [0,t_0]} |X^{(1)}(s)|^2\Big\}\le 4 \E \sup_{s\in [0,t_0]}  |X^{(1)}(s)|^2,$$ we obtain \eqref{*3}.
\end{proof}

\beg{proof}[Proof of Theorem \ref{T2.1}] Without loss of generality, we only consider  $s=0$ and simply denote $X_{0,}=X$; i.e.\,$X_{0,}(t)=X(t), X_{0,t}=X_t, t\ge 0$.

(1) Since the uniqueness follows from Theorem \ref{T2.1}(2),  which will be proved in the next step, in this step we only  prove   existence and estimate  \eqref{ES1}.
By Lemma \ref{L2.1},  there exists a unique adapted continuous process $(X_t)_{t\in [0,t_0]}$ such that
\beq\label{A01}\lim_{n\to\infty} \sup_{t\in [0,t_0]} \W_2(\mu_t^{(n)},\mu_t)^2\le \lim_{n\to\infty} \E \sup_{t\in [0,t_0]}  |X^{(n)}(t)- X(t)|^2=0,
\end{equation}
where $\mu_t$ is the distribution of $X_t$. By \eqref{EN},
$$X^{(n)}(t)=X(0)+ \int_0^t  b(s,X^{(n)}_s,\mu_s^{(n-1)}) \d s +\int_0^t\si(s,X^{(n)}_s,\mu_s^{(n-1)})\d W(s). $$ Then   \eqref{A01}, $(H1)$, $(H3)$  and the dominated convergence theorem imply that $\P$-a.s.
$$X(t)= X(0)+\int_0^t b(s,X_s, \mu_s)\d s +\int_0^t \si(s,X_s, \mu_s)\d W(s),\ \ t\in [0,t_0].$$
Therefore, $(X_t)_{t\in [0,t_0]}$ solves \eqref{ED1} up to time $t_0$, and \eqref{A01} implies  $ \E \sup_{s\in [0,t_0]} |X(s)|^2<\infty.$ The same holds for
$(X_{s,t})_{t\in [s,(s+t_0)\land T]}$ and $s\in [0,T]$. So, by solving the equation piecewise in time, and using the arbitrariness of $T>0$,
  we conclude that  \eqref{ED1} has a strong solution
  $(X_t)_{t\ge 0}$ with
  \beq\label{*X} \E \sup_{s\in [0,t]}|X(s)|^2<\infty,\ \ \ t\geq 0.\end{equation}

(2) By It\^o's formula and $(H2)$, we have
\begin{equation*}\begin{split}
\d \{\e^{\kk t}|X(t)-Y(t)|^2\}\le& 2\e^{\kk t} \big\<X(t)-Y(t), \{\si(t,X_t,\L_{X_t})-\si(t,Y_t,\L_{Y_t})\}\d W(t)\big\>\\
+&\e^{\kk t} \big\{\bb_1(t) \|X_t-Y_t\|_{\infty}^2+\bb_2(t)\W_2(\L_{X_t},\L_{Y_t})^2\big\}\d t.
\end{split}\end{equation*}  Noting that $\W_2(\L_{X_t},\L_{Y_t})^2\le \E\|X_t-Y_t\|_{\infty}^2$, we see that
$\gg_t:= \sup_{s\in [-r_0,t]} \e^{\kk s^+}|X(s)-Y(s)|^2$ satisfies
\beq\label{KD1} \beg{split} \E\gg_t \le &\E\|X_0-Y_0\|_\infty^2 + \E\int_0^t (\bb_1+\bb_2)(r)\e^{\kk r} \|X_r-Y_r\|_\infty^2\d s\\
&+2\E \sup_{s\in [0,t]} \int_0^s \e^{\kk r}  \big\<X(r)-Y(r), \{\si(r,X_r,\L_{X_r})-\si(r,Y_r,\L_{Y_r})\}\d W(r)\big\>.\end{split}\end{equation} By $(H1)$, the BDG inequality for $p = 1$ and since
 $\W_2(\L_{X_r},\L_{Y_r})^2\le \E\|X_r-Y_r\|_{\infty}^2$,  we have
\beg{align*} &   2\E \sup_{s\in [0,t]} \int_0^s \e^{\kk r}  \big\<X(r)-Y(r), \{\si(r,X_r,\L_{X_r})-\si(r,Y_r,\L_{Y_r})\}\d W(r)\big\>\\
&\le 4 \E\bigg(\int_0^t \e^{2\kk s} |X(s)-Y(s)|^2 \big(\aa_1(s)\|X_s-Y_s\|_\infty^2 +
\aa_s(s) W_2(\L_{X_r},\L_{Y_r})^2\big)\d s\bigg)^{\ff 1 2}\\
&\le \vv \E \gg_t + \ff 4 \vv \int_0^t (\aa_1(s)+\aa_2(s)) \E[\e^{\kk s}\|X_s-Y_s\|_\infty^2]\d s\\
&\le \vv \E\gg_t+ \ff 4 \vv \e^{\kk r_0} \int_0^t (\aa_1(s)+\aa_2(s)) \E\gg_s \d s.\end{align*}
Combining this with \eqref{KD1} we obtain
$$\E\gg_t\le \ff {\E\|X_0-Y_0\|_\infty^2}{1-\vv} + \ff{\e^{\kk r_0}}{1-\vv} \int_0^t \Big\{\ff{4}\vv (\aa_1(s)+\aa_2(s)) + \bb_1(s)+\bb_2(s)\Big\}\E\gg_s \d s,\ \ t\ge s.$$ So, Gronwall's Lemma implies
$$\E\gg_t\le \ff {\E\|X_0-Y_0\|_\infty^2}{1-\vv} \exp\bigg[\ff{\e^{\kk r_0}}{1-\vv} \int_0^t \Big\{\ff{4}\vv (\aa_1(s)+\aa_2(s)) + \bb_1(s)+\bb_2(s)\Big\} \d s\bigg].$$ Noting that $\E\gg_t\ge \e^{(t-r_0)\kk} \E\|X_t-Y_t\|_\infty^2$, this implies
 \beg{align*} &\E \|X_t-Y_t\|_\infty^2 \le \ff{\E \|X_0-Y_0\|_\infty^2}{1-\vv} \\
 &\times \exp\bigg[(r_0-t)\kk +\ff{\e^{\kk r_0}}{1-\vv} \int_0^t \Big\{\ff{4}\vv (\aa_1(s)+\aa_2(s)) + \bb_1(s)+\bb_2(s)\Big\} \d s\bigg].\end{align*}
 Since $(H2)$ remains true if $\kk$ is replaced by a smaller constant $\dd$, this estimate also holds for $\dd\in [0,\kk]$ replacing $\kk$. Therefore, the estimate in Theorem \ref{T2.1}(2) holds.

(3) Let   $(X_t)_{t\ge 0}$ solve \eqref{ED1} with $\scr L_{X_0}=\mu_0$, and let   $(\tt X_t,\tt W(t))$ on
$(\tt\OO, \{\tt\F_t\}_{t\ge 0}, \tt\P)$ be a weak solution of \eqref{ED1}  such that  $\L_{X_0}|_{\P}= \L_{\tt X_0}|_{\tt\P}=\mu_0$, i.e.
$\tt X_t$ solves
\beq\label{E1'} \d \tt X(t) = b(t,\tt X_t, \L_{\tt X_t}|_{\tt\P})\d t + \si(t,\tt X_t, \L_{\tt X_t}|_{\tt\P})\d \tt W(t),\ \ \ \scr L_{\tt X_0}=\mu_0.\end{equation}
 We aim to   prove   $\L_{X}|_{\P}=\L_{\tt X}|_{\tt\P}$.
  Let $\mu_t= \L_{X_t}|_{\P}$ and
$$\bar b(t,\xi)= b(t,\xi, \mu_t),\ \ \bar \si(t,\xi)= \si(t,\xi,\mu_t),\ \ t\ge 0, \xi\in\C.$$
  By $(H1)$-$(H3)$, the   stochastic functional differential equation
\beq\label{E10} \d \bar X(t) = \bar b(t,\bar X_t)\d t + \bar \si(t,\bar X_t)\d \tt W(t),\ \ \bar X_0= \tt X_0 \end{equation}
has a unique solution.
According to Yamada--Watanabe, it also satisfies weak uniqueness. Noting that
$$\d X(t) = \bar b(t, X_t)\d t + \bar \si(t,  X_t)\d  W(t),\ \ \L_{X_0}|_{\P}= \L_{\tt X_0}|_{\tt\P},$$ the weak uniqueness of \eqref{E10} implies
\beq\label{HW} \L_{\bar X}|_{\tt\P}= \L_X|_{\P}.\end{equation}
So, \eqref{E10} reduces to
$$ \d \bar X(t) = b(t,\bar X_t, \L_{\bar X_t}|_{\tt\P})\d t + \si(t,\bar X_t, \L_{\bar X_t}|_{\tt\P})\d \tt W(t),\ \ \bar X_0=\tt X_0.$$
Since the strong uniqueness of   \eqref{E1'} is ensured by Step (1),   we obtain  $\bar X=\tt X$. Therefore,  \eqref{HW} implies $\L_{\tt X}|_{\tt \P} = \L_X|_{\P}$ as wanted.

Finally, since $\C$ is a Polish space, for any $\mu_0,\nu_0\in \scr P_2^\C$, we can take $\F_0$-measurable random variables $X_0$, $Y_0$ such that $\L_{X_0}=\mu_0, \L_{Y_0}=\nu_0$ and
$\W_2(\mu_0, \nu_0)^2=\E\|X_0-Y_0\|_{\infty}^2$. Combining this with  $\W_2(P_t^*\mu_0, P_t^*\nu_0)^2\le \E\|X_t-Y_t\|_{\infty}^2,$  we deduce the estimate in Theorem \ref{T2.1}(3) from
that in Theorem \ref{T2.1}(2).
\end{proof}

\section{Harnack inequality and applications}

To prove Theorem \ref{T1.2}, we investigate   Harnack inequalities of the operator $P_t$ defined by
\beq\label{*D2} (P_tf)(\mu_0)= \int_{\C} f  \d(P_t^*\mu_0),\ \ f\in \B_b(\C), t\ge 0, \mu_0\in \scr P_2^\C.\end{equation}
We will consider the Harnack inequality with a power $p>1$ introduced in \cite{W97}, and the log-Harnack inequality developed in \cite{RW10, W10}, where classical SDEs on $\R^d$ and manifolds are considered.
To establish these inequalities for the present path-distribution dependent SDEs, we will adopt coupling by change of measures
introduced in \cite{ATW06,W07}. We refer to \cite{Wbook} for a general theory on this method and applications.

To construct the desired coupling for the segment solution $X_t$, we need to assume that $\si(t,\xi,\mu)=\si(t,\xi(0))$; that is, we consider the following simpler version of \eqref{ED1}:
\begin{equation}\label{E11}
\d X(t)= b(t,X_t,\scr L_{X_t})\d t +\si(t,X(t))\d W(t).
\end{equation}

\beg{thm}\label{T3.1} Assume {\bf (A)}. Then there exists $H_1\in C(\R_+;\R_+)$ such that
  for any    $ \mu_0,\nu_0\in \scr P_2^\C $,    $\F_0$-measurable random variables
$X_0, Y_0$ with $\L_{X_0}=\mu_0, \L_{Y_0}=\nu_0$,   and   $f\in \B_{b}^+(\C)$,
 \beq\label{LH11}
(P_{T}\log f)(\nu_0)\le  \log (P_{T}f)(\mu_0)+ H_1(T) \mathbb{E}\bigg( \ff{|X(0)-Y(0)|^2}{T-r_0} + \|X_0-Y_0\|_\infty^2 \bigg), \ \ T>r_0.
\end{equation}
 If moreover $\eqref{*P}$ holds for some increasing $\kk_3:\R_+\to\R_+$, then there exists $H_2\in C(D;\R_+)$, where $D$ is as in Theorem \ref{T1.2},  such that
 \beq\label{HI11}   (P_{T}f)(\nu_0)\le (P_{T}f^p)^{\frac{1}{p}}(\mu_0)\mathbb{E}\Big(\e^{H_2(p,T) \big(1+\ff{|X(0)-Y(0)|^2}{T-r_0} + \|X_0-Y_0\|_\infty^2\big)}\Big),\ \ T>r_0, (p,T)\in D
   \end{equation}
holds for   $\mu_0,\nu_0$ and $X_0,Y_0$ as above.
\end{thm}

As a consequence of Theorem \ref{T3.1}, we have the following result, see, for instance,  the proof of \cite[Prposition 3.1]{WY11}.

\beg{cor}\label{C3.3} Assume {\bf (A)} and let $T> r_0$. For any $\mu_0,\nu_0\in\scr P_2$, $P_{T}^*\mu_0$ and $P_{T}^*\nu_0$ are equivalent and the Radon-Nykodim derivative satisfies
the entropy estimate
$$
\int_{\C} \bigg(\log \ff{\d P_{T}^*\nu_0}{\d P_{T}^*\mu_0}\bigg)\d P_{T}^*\nu_0\le\inf_{\L_{X_0}=\mu_0, \scr L_{Y_0}=\nu_0} \mathbb{E}\bigg[H_1(T)\bigg(\ff{|X(0)-Y(0)|^2}{T-r_0} +\|X_0-Y_0\|_\infty^2\bigg)\bigg] ,\ \ T>r_0. $$
If $\eqref{*P}$ holds, then for any $T>r_0$ and $p>(1+\kk_3(T)\ll(T))^2$, $$\int_{\C}  \bigg(\ff{\d P_{T}^*\nu_0}{\d P_{T}^*\mu_0}\bigg)^{\ff 1 p} \d (P_{T}^*\nu_0)
 \le \inf_{\L_{X_0}=\mu_0, \scr L_{Y_0}=\nu_0}  \mathbb{E}\Big(\e^{H_2(p,T) \big(1+\ff{|X(0)-Y(0)|^2}{T-r_0} + \|X_0-Y_0\|_\infty^2\big)}\Big).$$
\end{cor}

\

\beg{proof}[Proof of Theorem \ref{T3.1}] For $\mu_t:= P_t^*\mu_0$ and $\nu_t:=P_t^*\nu_0$, we may
rewrite \eqref{E11} as
\begin{equation}\label{barX}
\d X(t)= \bar{b}(t,X_t)\d t +\si(t,X(t))\d \bar{W}(t),\ \ \scr L_{X_0}=\mu_0,
\end{equation}
where
\begin{equation*}\begin{split}
&\bar{b}(t,\xi):=b(t,\xi,\nu_t),\ \ \d \bar{W}(t):=\d W(t)+ \bar{\gamma}(t)\d t,\\
& \bar{\gamma}(t):=\si^{-1}(t,X(t))[b(t,X_t,\mu_t)- b(t,X_t,\nu_t)].
\end{split}\end{equation*}
By assumption {\bf (A)}  and Theorem \ref{T2.1}(3), we have
\begin{equation}\begin{split}\label{EbarG}
|\bar{\gamma}(t)|\leq \lambda(t) \kk_2(t)\W_2(\mu_t,\nu_t)\le K(t) \W_2(\mu_0,\nu_0), \ \ t\in[0,T]
\end{split}\end{equation} for some increasing function $K:\R_+\to\R_+.$
Let
\begin{equation}\label{EB2}
\bar{R}_t=\exp\left\{-\int_0^t\langle\bar{\gamma}(s),\d W(s)\rangle-\frac{1}{2}\int_0^t|\bar{\gamma}(s)|^2\d s\right\},\ \ t\in[0,T].
 \end{equation}
  By Girsanov's theorem, $\{\bar{W}(t)\}_{t\in[0,T]}$ is a $d$-dimensional Brownian motion under the probability measure $\bar{\P}_T:= \bar R_T\P$.

Next, according to the proof of \cite[Theorem 4.3.1]{Wbook} or \cite[Theorem 1.1]{WY11}, we can construct an adapted process $\tt\gg(t)$ on $\R^d$ such that
\beg{enumerate} \item[(a)] Under the probability measure $\bar \P_T$,
$$\tilde{R}_t:=\exp\left\{-\int_0^t\langle\tilde{\gamma}(s),\d \bar{W}(s)\rangle-\frac{1}{2}\int_0^t|\tilde{\gamma}(s)|^2\d s\right\}, \ \ t\in[0,T]$$
is a martingale, such that $\tt\P_T:= \tt R_T\bar\P_T=\tt R_T\bar R_T\P$ is a probability measure under which
$$\tt W(t):=  \bar{W}(t)+ \int_0^t \tilde{\gamma}(s)\d s =  W(t)+ \int_0^t \big(\bar{\gamma}(s)+\tilde{\gamma}(s)\big)\d s,\ \ t\in [0,T]$$ is a $d$-dimensional Brownian motion.
\item[(b)] Letting $Y(t)$ solve the following stochastic functional differential equation under the probability measure $\tt\P_T$ with the given initial value $Y_0$:
\beq\label{barCY}\begin{split}
\d Y(t) &= \bar{b}(t,Y_t)\d t+ \si(t,Y(t)) \d \tilde{W}(t),
\end{split}\end{equation}  we have  $\L_{Y_0|\tt\P}=\L_{Y_0}=\nu_0$ and  $X_T=Y_T\ \tt\P_T$-a.s.
\item[(c)] There exists  $C\in C(\R_+;\R_+)$ such that
$$ \E_{\tt\P_T} \int_0^T |\tt\gg(s)|^2\d s \le C(T) \E\Big(\ff{|X(0)-Y(0)|^2}{T-r_0}+\|X_0-Y_0\|_\infty^2\Big).$$
    \end{enumerate}

By the definition of $\bar b$ we see that $(Y_t,\tt W(t))$ is a weak solution to the equation \eqref{barX} with initial distribution $\nu_0$, so that by the weak uniqueness,
$\scr L_{Y_t}|_{\tt \P_T}=\nu_t, t\in [0,T].$ Combining this with (b) we obtain
$$(P_Tf)(\nu_0)= \E_{\tt \P_T} [f(Y_T)] = \E_{\tt P_T}[f(X_T)] = \E [\bar R_T\tt R_T f(X_T)],\ \ f\in \B_b^+(\C).$$
Letting $R_T= \bar R_T\tt R_T$,   by Young's inequality and H\"older's inequality respectively, we obtain
\beq\label{LHI}  (P_T\log f)(\nu_0) \le \E [  R_T \log   R_T ]+ \log\E[f(X_T)]
 = \E [  R_T \log   R_T ]+ \log(P_T f)(\mu_0),
 \end{equation}
and
\beq\label{HI}\begin{split}
&(P_T f(\nu_0))^{p}\le (\E R_T^{\ff {p}{p-1}})^{p-1} (\E f^{p}(X_T))= (\E R_T^{\ff {p}{p-1}})^{p-1} P_T  f^{p}(\mu_0),\ \ p>1.
\end{split}\end{equation}
 We are now ready to prove assertions (1) and (2) as follows.

 By \eqref{EbarG} ,  (c)   and since $\W_2(\mu_0,\nu_0)^2\le \E\|X_0-Y_0\|_\infty^2$,
\begin{equation}\label{LR}\beg{split}
\E[R_T\log R_T]&\le \frac{1}{2}\E_{\tilde{\P}_T}\int_0^T|\bar{\gamma}(s)+\tilde{\gamma}(s)|^2\d s\\
&\le \E_{\tilde{\P}_T}\int_0^T|\tilde{\gamma}(s)|^2\d s+ \int_0^T|\bar{\gamma}(s)|^2\d s\\
&\le \E_{\tilde{\P}_T}\int_0^T|\tilde{\gamma}(s)|^2\d s+\int_0^T\lambda(t)^2\kk_2(t)^2\W_2(\mu_t,\nu_t)^2\d t\\
&\le H_1(T) \mathbb{E}\bigg(\ff{|X(0)-Y(0)|^2}{T-r_0}+ \|X_0-Y_0\|_\infty^2\bigg), \ \ T>r_0\end{split}\end{equation} holds for some $H_1\in C(\R_+;\R_+).$
Combining this with \eqref{LHI} we obtain \eqref{LH11}.

 Finally, according to the proof of \cite[Theorem 4.1]{WY11}, there exists $C\in C(D;\R_+)$ such that
$$(\E_{\bar\P_T} \tt R_T^{\ff {p}{p-1}})^{\frac{p-1}{p}} \le   \E  \Big(\e^{C(p,T)\big(1+\ff{|X(0-Y(0)|^2}{T-r_0} +\|X_0-Y_0\|_\infty^2\big)}\Big),\ \ T>r_0, (p,T)\in D.  $$
For any $p>p(T):= (1+\kk_3(T)\ll(T))^2$, by applying this estimate for $p_1:=\ff 1 2(p+(p(T))$ and combining with $R_T=\tt R_T \bar R_T$, \eqref{EbarG}, \eqref{EB2} and
$\W_2(\mu_0,\nu_0)^2\le \E\|X_0-Y_0\|_\infty^2,$ we arrive at
\beg{align*} \Big(\E  R_T^{\ff {p}{p-1}}\Big)^{\frac{p-1}{p}} &=\Big(\E_{\bar\P_T}   \tt R_T^{\ff {p}{p-1}} \bar R_T^{\ff 1 {p-1}}\Big)^{\frac{p-1}{p}}
 \le  \Big(\E_{\bar \P_T} \tt R_T^{\ff {p_1}{p_1-1}}\Big)^{\frac{p_1-1}{p_1}} \Big(\E_{\bar\P_T} \bar R_T^{\ff {p_1}{p-p_1}}\Big)^{\frac{p-p_1}{pp_1}}\\
&\le  \E  \Big(\e^{C(p_1,T)\big(1+\ff{|X(0-Y(0)|^2}{T-r_0} +\|X_0-Y_0\|_\infty^2\big)}\Big) \Big(\E\bar R_T^{\ff {p}{p-p_1}}\Big)^{\frac{p-p_1}{pp_1}}\\
&\le \E \Big(\e^{H_2(p,T)\big(1+\ff{|X(0-Y(0)|^2}{T-r_0} +\|X_0-Y_0\|_\infty^2\big)}\Big),\ \ T>r_0, (p,T)\in D\end{align*} for some $H_2\in C(\R_+;\R_+).$ Therefore, \eqref{HI11} follows from \eqref{HI}.
\end{proof}

\section{Shift Harnack inequality and integration by parts formula}

To prove Theorem \ref{T1.3}, we investigate the shift Harnack inequality and integration by parts formula introduced  in \cite{W14a}.     Assume that $\si(t,\xi,\mu)=\si(t)$ is
invertible. Then the path-distribution dependent SDE \eqref{ED1} becomes
$$
\d X(t)= b(t,X_t, \L_{X_t})\d t +\si(t) \d W(t),\ \ \scr L_{X_0}=\mu_0.$$
 To apply the existing shift Harnack inequality and integration by parts formula, we let
 $$\bar b(t,\xi):= b(t, \xi,\mu_t),\ \ \mu_t:=\scr L_{X_t}=P_t^*\mu_0,\ \ t\ge 0,\xi\in \C$$ and
 rewrite this equation as
$$\d X(t)= \bar b(t,X_t) \d t +\si(t) \d W(t),\ \ \ \scr L_{X_0}=\mu_0.$$
Then the following result follows from   \cite[Theorem 4.2.3]{Wbook}.

\beg{thm}\label{T4.1} Let $\si: [0,\infty)\to \R^d\otimes \R^d$ and $b: [0,\infty)\times \C\times\scr P_2^\C\to \R^d$ satisfy {\bf (A)}, and assume that for any
$(t,\mu)\in \R_+\times \scr P_2^\C$, $b(t,\cdot,\mu)$ is differentiable. Then
$$\LL(T):=\sup_{t\in[0,T]}\|\si(t)^{-1}\|^2<\infty,\ \ K(T):= \sup_{t\in [0,T],\mu\in \scr P_2^\C} \|\nn b(t,\cdot,\mu)\|_\infty^2<\infty,\ \   T\geq 0.$$  Moreover:
\beg{enumerate} \item[$(1)$] For any $p>1,T>r_0, \mu_0\in \scr P_2^\C, \eta\in\mathbb{H}^{1} $ and $f\in \B_b^+(\C)$,
\beg{align*}
(P_{T}f)^p(\mu_0)\le &(P_{T}f^p(\eta+\cdot))(\mu_0)\\
&\times \exp\bigg[\ff{p\, \LL(T)\left(1
+T^2K(T) \right)\left(\frac{|\eta(-r_0)|^2}{T-r_0}+\|\eta\|_{\mathbb{H}^{1}}^{2}\right) }{(p-1)^2}\bigg],\ \ p>1,
\end{align*}
and
$$
(P_{T}\log f)(\mu_0)\le \log (P_{T} f(\eta+\cdot))(\mu_0)
+ \LL(T)\left(1
+T^2K(T) \right)\left(\frac{|\eta(-r_0)|^2}{T-r_0}+\|\eta\|_{\mathbb{H}^{1}}^{2}\right).
$$
\item[$(2)$] For any $T>r_0,$ let   \begin{equation*}\begin{split}
&\Phi (t)=1_{[0,T-r_0]}(t)\ff{\eta(-r_{0})}{T-r_0}+1_{(T-r_0,T]}(t)\eta^{'}(t-T),\\ &\Theta (t)=\int_0^{t^+}\Phi (s)\d s,\ \ t\in[-r_0,T].\end{split}\end{equation*}
Then for any  $f\in C^1(\C)$,   $\eta\in \mathbb{H}^{1}$ and $\mu_0\in \scr P_2^\C$,
$$\E(\nn_\eta f)(X_{T})= \E\bigg[f(X_{T})\int_0^T  \big\<\si(t)^{-1}(\Phi (t) -\nn_{\Theta_t} b(t,\cdot, P_{t}^*\mu_0)(X_t)),\ \d W(t)\big\>\bigg].$$
\end{enumerate} \end{thm}
 As consequence of Theorem \ref{T4.1} we have the following result.

\beg{cor}\label{C4.2} In the situation of Theorem \ref{T4.1}. For any $\mu_0\in\scr P_2^\C, \eta\in\H^1$ and $T>r_0$, $\mu_T:= P_T^*\mu_0$ satisfies
\beg{align*}&\int_\C \Big(\log \ff{\d \mu_T(\cdot+\eta)}{\d\mu_T}\Big)\d\mu_T(\cdot+\eta)\le \LL(T) (1+T^2 K(T))\Big(\ff{|\eta(-r_0)|^2}{T-r_0}+ \|\eta\|_{\H^1}^2\Big),\\
&\int_\C \Big( \ff{\d \mu_T(\cdot+\eta)}{\d\mu_T}\Big)^{\ff 1 p} \d\mu_T(\cdot+\eta)\le \exp\bigg[\ff{\LL(T) (1+T^2 K(T))}{(p-1)^2}\Big(\ff{|\eta(-r_0)|^2}{T-r_0}+ \|\eta\|_{\H^1}^2\Big) \bigg],\ \ p>1,\\
&\int_\C \Big|\ff{\d \pp_\eta\mu_T}{\d\mu_T}\Big|^2\d\mu_T \le \LL(T) \big(1+K(T)T^2\big) \Big(\ff{|\eta(-r_0)|^2}{T-r_0}+ \|\eta\|_{\H^1}^2\Big).\end{align*}
\end{cor}

\beg{proof}
The first two estimates follow from  Theorem \ref{T4.1}(1), see    \cite{W14a} or \cite[\S 1.4]{Wbook}. As the last estimate is not explicitly given in these references, we present a brief proof below.
It is easy to see that $$M(T):= \int_0^T  \big\<\si(t)^{-1}(\Phi (t) -\nn_{\Theta_t} b(t,\cdot, P_{t}^*\mu_0)(X_t)),\ \d W(t)\big\>$$ satisfies
$$\E M(T)^2\le C(T):= \LL(T) \big(1+K(T)T^2\big) \Big(\ff{|\eta(-r_0)|^2}{T-r_0}+ \|\eta\|_{\H^1}^2\Big).$$
Then, Theorem \ref{T4.1}(2) implies that
$$C^1(\C)\ni f\mapsto (\pp_\eta \mu_T)(f):= \bigg(\ff{\d }{\d\vv} \int_\C f\d\mu_T(\cdot+\vv \eta)\bigg)\bigg|_{\vv=0}$$ is a densely defined bounded linear functional on $L^2(\mu_T)$ with
$$\big|(\pp_\eta \mu_T)(f)\big|^2 \le \mu_T(f^2)\E M(T)^2\le C(T)\mu_T(f^2).$$
By the Riesz Representation Theorem, it uniquely extends to a bounded linear functional
$$(\pp_\eta\mu_T)(f):= \int_\C fg\d\mu_T,\ \ \ f\in L^2(\mu_T)$$  for some $g\in L^2(\mu_T)$ with $\mu_T(g^2)\le C(T).$ Consequently,
$\mu_T$ is differentiable along $\eta$ with  $(\pp_\eta\mu_T)(A)=\int_Ag\d\mu_T, A\in \B(\C),$ and $\pp_\eta \mu_T$ is absolutely continuous with respect to $\mu_T$ such that
$$\int_\C \Big(\ff{\d\pp_\eta\mu_T}{\d\mu_T}\Big)^2\d\mu_T= \int_\C g^2\d\mu_T\le C(T).$$\end{proof}

\paragraph{Acknowledgement.} The authors would like to thank the referee for helpful comments on an earlier version of the paper.

\end{document}